\documentclass[11pt]{article}

\usepackage[T2A]{fontenc}
\usepackage[cp1251]{inputenc}
\usepackage[tbtags]{amsmath}
\usepackage{amsfonts,amsthm,amssymb,amstext,amsbsy,
latexsym,mathrsfs,mparhack,marginnote}

\usepackage{mathrsfs}

\usepackage{graphicx}

\usepackage{accents}

\usepackage[dvipsnames]{xcolor} 

\usepackage[bbgreekl]{mathbbol}
\usepackage{bold-extra}

\DeclareSymbolFontAlphabet{\mathbbl}{bbold}

\usepackage{hyperref}
\hypersetup{citecolor=magenta}

\input{00.sty}

\begin{document}

\title
{Locally countable graphs 
of second projective class 
not generated by countably many projective 
functions %
\vyk
{The research was carried out within the state 
assignment of Ministry of Science and Higher 
Education of the Russian Federation for IITP RAS.}
}

\author 
{
Vladimir~Kanovei\thanks{IITP RAS,
  Moscow, Russia, \ 
  {\tt kanovei@iitp.ru} --- contact author. 
}  
\and
Vassily~Lyubetsky\thanks{IITP RAS,
  Moscow, Russia, \ {\tt lyubetsk@iitp.ru} 
}
}

\date 
{\today}

\maketitle

\markboth{В.\,Г.~Кановей и В.\,А.~Любецкий}
{Locally countable graphs  
of 2nd projective class not generated 
by projective functions}

\begin{abstract}
To answer a question by Rettich and Serafin, 
we define a model of set theory in which there 
exists a locally countable $\ip12$ graph  
on the reals, which is not 
generated by a countable family of projective 
(or even real-ordinal definable, ROD) functions.
We also prove that the $\is12$ 
equi-constructibility graph on the reals 
is not generated by a countable family of ROD 
functions in the Solovay model.
\end{abstract}

\vyk{
}

\parf{Introduction}
\las{int}

A graph $G$ is said to be 
\rit{generated by a family of functions} $F$, 
$G=G_F$ for brevity,   
\cite[\S\,2]{chrom} or \cite{rettich}, 
if for any $x,y$ in the domain of $G$, it holds 
{%
\bul
{e*}
{%
{x \mathrel G y}
\iff 
\sus f\in F\big(x=f(y) \lor y=f(x)\big).
}%
\setcounter{saved}{\value{equation}}%
}%
Such a graph has to be locally countable, of course, 
provided the family $F$ is countable or finite and 
each $f\in F$ is ${\le}\alo$-to-1 
(\ie, the $f$-preimage of every element is at most countable). 
Conversely, using the axiom of choice, 
one easily proves that every locally countable 
graph is generated by a countable 
family of functions, by necessity of type 
${\le}\alo$-to-1. 

The inverse problem becomes more involved in
the case where definability
of generating functions is required depending 
on the definability of a given graph. 
In particular, as shown in \cite[\S\,2]{chrom}, 
every locally countable Borel graph is generated 
by a countable family of Borel functions. 
(See Lemma~\ref{pro1} below.)

\vyk{
This is because any  Borel set with countable 
sections is a countable union of Borel 
uniform sets by a well-known theorem of classical 
descriptive set theory. 
}

This paper is devoted to the problem of generation 
of graphs of second projective level by 
countable systems 
of real-ordinal definable 
(ROD, for brevity) functions, 
including projective  functions of any level. 
Both of our main results 
(Theorems \ref{mts} and \ref{mt} below)
state that such generation is impossible in certain
situations. 
A simple reduction (Lemma~\ref{22} below) allows
us to also deduce the nonexistence of 
a countable ROD
coloring in those cases.

At the $\fs12$ side, we prove 
(Theorem~\ref{t2}) that every locally 
countable $\fs12$ graph is generated 
by a family of $\ali$-many Borel functions, 
by necessity of type ${\le}\alo$-to-1. 
Then we consider the equi-constructibility 
$\is12$ graph 
{\setcounter{equation}{\value{saved}}
\bul
{e|}
{%
{x\mcL y} \ \ := \ \ {x\ne y\land \rL[x]=\rL[y]}
\quad-\quad \text{for }\, x,y\in\dn;
}%
}%
the irreflexive part of the equivalence 
relation of equi-constructibility on $\dn.$ 
Its properties depend on a model in which $\cL$ 
is considered. 
(See \eg\ \cite{payne,rinot17} on different properties 
of one and the same graph in different set universes.)
In particular, in $\rL$, the constructible 
universe, any reals $x\ne y$ in $\dn$ 
are $\cL$-adjacent, and, by Theorem~\ref{t2+} 
below, 
$\cL$ is generated by a countable family of 
$\is12$ functions 
(not necessarily of type ${\le}\alo$-to-1). 


A similar positive result also holds for some 
models of the form $\rL[a]$, for instance, in 
the case of Sacks-generic reals $a\in\dn.$ 
However Lemma~\ref{sax} solves the issue in the 
negative provided we require 
ROD countable families of  
functions to generate $\cL$. 
The lengthy proof of the lemma involves some 
specific properties of the Sacks forcing, and 
the result itself is the key ingredient in the 
proof of the following theorem. 
It involves  \rit{the Solovay model}, \ie,  
a model of the theory $\zfc$ 
+ ``all ROD sets
of reals are Lebesgue measurable'',
defined in \cite{sol}.
See \S\,\ref{solm} for more information 
about this model.

\bte
\lam{mts}
It is true in the Solovay model that\/ $\cL$ 
is a locally countable\/ $\is12$ graph  
not generated by countably many\/ {\em ROD} 
functions, and\/ $\cL$ has the 
$\ROD$-chromatic number  
$\chi_{\text{\sc rod}}(\cL)=\gc=2^{\alo}.$ 
\ete

As usual in descriptive set theory, the case 
of $\fp12$ graphs causes even greater difficulties. 
The following theorem is our second 
main result. 

\bte
\lam{mt}
There exists a model of set theory\/ $\zfc$ 
in which there is a 
locally countable\/ $\ip12$ 
graph\/ $G$ on\/ $\dn$  
not generated by countably many\/ $\ROD$  
functions, and\/ $G$ has the 
$\ROD$-chromatic number\/  
$\chi_{\text{\sc rod}}(G)=\gc=2^{\alo}.$  
\ete

\vyk{
We recall that countable equivalence relations 
are those which have at most countable 
equivalence classes (see \eg\ \cite{K2024}). 
This corresponds to locally countable   
graphs, of course. 
}

Theorem \ref{mt} answers in the positive 
the question on the existence of such a model 
and a graph in \cite[Problem 3.6]{rettich}, 
with the required example of a graph presented 
in the special form of 
the irreflexive part of a countable $\ip12$ 
equivalence relation. 
The proof of this result in \S\,\ref{p12} 
makes use of a model defined in \cite{kl28}, 
in which there exists a planar $\ip12$ set 
with countable sections, not uniformizable by 
ROD (including projective) sets. 

The ``moreover'' parts of the theorems deal 
with the $\ROD$-chromatic numbers 
(\ie\ those emerging from $\ROD$ colorings)
of the graphs considered, in both cases 
equal to the continuum 
$\gc=2^{\alo},$ while the classical
chromatic numbers do not exceed $\alo$ due to 
their local countability. 
The results on the separation of classical and Borel
chromatic numbers of various graphs 
are known from a number of works, see \eg\  
\cite{bern,dPT,chrom, millBgraph,rettich} 
among other papers, 
and also an unpublished book by Kechris and 
Marks on descriptive graph combinatorics\snos
{A.\,S.\;Kechris and A.\,S.\;Marks,
{\em Descriptive graph combinatorics\/}, 
Preliminary version, 2020,
\url{https://math.berkeley.edu/~marks/papers/combinatorics20book.pdf}
  (accessed 27.May.2026).}. 
Among the more special directions in 
this research domain, a very recent paper 
\cite{gre2026} develops a pretty 
novel technique of using forcing to define Borel 
and analytic {\em directed\/} 
graphs with various effects 
related to Borel chromatic numbers. 
Another direction, presented in \cite{marks},
is related to applications of the technique 
of infinite games (determinateness). 
\rit{Lebesgue measurable} colorings and
corresponding chromatic numbers are considered
in comparison with classical  and Borel 
colorings for some graphs $G$ in  
\cite{marks23,payne,rinot12,rinot17,sosh1} and others.

In this overall picture, 
оur results open up a new direction in 
descriptive graph combinatorics, and give
some stronger results, since Borel
maps make up only a small part of $\ROD$ mappings. 
For the importance of research in this area, see,
for example, \cite[\S\S\,4 and 5]{rettich}.

\parf{Preliminaries} 
\las{prel}

The reader is assumed to have a basic 
knowledge of descriptive set theory and forcing. 
However we have to review here some 
more special notions.\vom 

\rit{Descriptive set theory.} 
As is customary, we take 
\rit{the Baire space} $\bai$ 
as the principal domain, and also consider 
$\bai\ti\bai$ \etc, as well as the 
\rit{Cantor discontinuum} $\dn\sq\bai.$ 
Elements of $\bai$ are called \rit{reals}.
Sets $X\sq\bai,$ $X\sq(\bai\ti\bai)$, \etc \ 
are \rit{pointsets}, collections of 
pointsets are \rit{pointclasses}. 
In particular, such pointclasses as 
$\fd11$ (Borel sets), 
$\fs11$ (analytic or Suslin sets), 
$\fp11$ (coanalytic or co-Suslin sets), 
$\fs1n\yi\fp1n\yi\fd1n$ (projective classes), 
$\fs1\iy=\fp1\iy$ (all projective sets) are 
systematically studied by descriptive 
set theory 
(see \eg\ \cite{Kdst} or 
\cite[Chap.\,11,\,25,\,32]{jechmill}).\vom

\rit{Ordinal definability.} 
Such pointclasses as 
$\OD$ (ordinal-definable sets), $\OD(P)$ 
(elements $p\in P$ are admitted as parameters 
along with ordinals), and 
$\ROD=\OD(\bai)$ (real-ordinal-definable sets) 
are considered in more abstract branches of set 
theory. 
See \cite[Chap.\,13]{jechmill} or \cite{sol} 
on OD and ROD.\vom

\rit{Tuples.} 
$\bse$ is the set of all (finite) \rit{tuples} 
of numbers $0$ and $1$, including the empty 
tuple $\etu$. 
The \rit{length} of $s\in\nse$ is denoted by 
$\lh s$. 
If $s\in\nse$ and $j=0,1$ then $s\we j\in\bse$ 
is obtained by adding $j$ as the new 
rightmost term. 
Generally, $s\we t$ is the concatenation. 
Then, $s\sq t$ means that a tuple $t$ 
extends $s$, 
and $\su$ means proper extension.\vom

\rit{Trees and perfect sets.} 
A set $T\sq\bse$ is a \rit{tree} iff 
(1) $s\su t\in T\imp s\in T$, and 
(2) $s\in T\imp \sus j\,(s\we j\in T)$, 
and a \rit{perfect} tree if additionally 
\bce
$\kaz s\in T\,\sus t\in T\,\sus j\ne i\in\ans{0,1}\,
(t\we i\in T\land t\we j\in T)$.
\ece
Then 
$[T]=\ens{x\in\dn}{\kaz m\,(x\res m\in T)}$ 
is a perfect subset of $\dn.$\vom 

\vyk{
In this case we let $\rat[T]$ 
(the ``rational points'' in $[T]$) 
be the set of all leftmost and rightmost 
branches in subsets of the form 
$\ens{x\in[T]}{s\su x}$, where $s\in T$. 
Then $\rat[T]\sq[T]$ is a countable dense subset 
of $[T]$.\vom

\rit{Coding continuous functions.}  
Assume that $T\sq\bse$ is a perfect tree and 
$f:[T]\to\bai$ is continuous. 
Then $f$ is fully determined by the restricted 
function $f\res \rat[T]$, which can be 
effectively coded by an appropriate real 
$c\in\bai$ because $\rat[T]\sq[T]$ is 
a countable set.\vom
}

\rit{Graphs.} 
A \rit{graph} will be any 
$G\sq\dn\ti\dn,$ 
\rit{symmetric} 
(${x\mre G y}\eqv {y\mre G x}$), 
and \rit{irreflexive} 
(${x\mre G y}\imp {x\ne y}$). 
Here and below ${x\mre G y}$ is a shortcut for 
$\ang{x,y}\in G$. 
A graph $G$ is \rit{of class $K$} 
(say $K=\fd11=$ Borel) if $G$ belongs 
to $K$ as a set of pairs. 
Thus a Borel graph $G$ means that $G$ is 
a Borel set.

Elements of the set 
$|G|=\ens{x}{\sus y\,(x\mre G y)}$ 
are \rit{vertices} of $G$ while the pairs 
in $G$ are its \rit{edges}. 
Occasionally, graphs $G$ with 
$|G|\sq{\bai\ti\bai}$ will be 
considered as well. 
A graph $G$ is \rit{locally countable} if for 
any vertice $x\in |G|$ the set 
$\ens{y}{x\mre G y}$ of all \rit{adjacent} 
vertices is at most countable. 
This is equivalent to the countability of every 
\rit{connected component} of $G$.\vom 

\rit{Generated graphs.} 
The notion of a graph $G$ 
\rit{generated by a family of functions} $F$, 
based on \eqref{e*}, was given in Introduction. 
Here we suppose that each $f\in F$ is a 
function $f:\dom f\sq|G|\to |G|$, and 
$\dom f=|G|$ is not assumed, while 
an equality like $y=f(x)$ means 
$x\in\dom f\land y=f(x)$.\vom 

\rit{$\ROD$ chromatic numbers.} 
The $\ROD$ chromatic number 
$\chi_{\text{\sc rod}}(G)$  
of a graph $G$ is 
the least possible cardinal of a set $\ran c$, 
where $c:|G|\to\dn$ is a $\ROD$ colouring of $G$. 
(Meaning that $c$ is a colouring in the classical 
sense and $c$ is a $\ROD$ map.)
Borel chromatic numbers $\chi_{\text{\sc bor}}(G)$ 
are defined similarly. 
Obviously
\bul
{chrn}
{
\chi(G)
\le
\chi_{\text{\sc rod}}(G)
\le
\chi_{\text{\sc bor}}(G)
\le\gc.
}
\vyk{
Another generalization option is connected with 
the \rit{Lebesgue measurable} colorings and
corresponding chromatic numbers, considered
in comparison with $\chi(G)$ and 
$\chi_{\text{\sc bor}}(G)$
for some graphs $G$ in  
\cite{marks23,payne,rinot12,rinot17,sosh1} and others.
}

\parf{Some preliminary results} 
\las{s12}

Here we present a few simple lemmas,
which (except for Lemma~\ref{22}) 
are not directly used in the proofs
of Theorems \ref{mts} and \ref{mt} 
(our main results), but contain some limitations
that allow to better evaluate the content
of our theorems. 

\ble
\lam{22}
Any\/ $\ROD$ graph\/ $G\sq\dn\ti\dn$ that allows
a countable\/ $\ROD$ coloring is generated by 
a countable family of\/ $\ROD$ functions.  

Similarly for Borel graphs,
countable Borel colorings, and countable
families of Borel functions.  
\ele

\bpf
Let $c:|G|\to\om $ be a $\ROD$-coloring.  
Let $x\in|G|$.  
The set
$G(x)=\ens{y}{x\mG y}$ 
all $G$-neighbors are at most countable, 
because $c$ is bijective on $G(x)$ by
the coloring definition. 
Therefore, the set 
$U_x=\ens{c(y)}{x\mG y}\sq\om$ has the form
$$
\bay{rcl}
U_x &=& \ans{m_0(x)<m_1(x)<m_2(x)<\dots},
\quad\text{or}\\[0.5ex]
U_x&=&
\ans{m_0(x)<m_1(x)<\ldots<m_k=m_{k+1}=m_{k+2}
=\dots}.
\eay
$$ 
This allows us to define $f_j(x)=y_j$, where 
$y_j\in G(x)$ is the only real 
satisfying $c(y_j)=m_j(x)$.
Then each $f_j:|G|\to|G|$ is a 
$\ROD$ map,
and it is easy to see that the graph $G$ is 
generated by the family $F=\ens{f_j}{j\ge1}$. 
\epf

The following two lemmas concern Borel graphs.

\ble
[{\cite[\S\,2]{chrom}}]
\lam{pro1}
Any locally countable Borel graph\/ $G$
is generated by a countable family\/  
$F$ of Borel functions. 
\ele

\bpf
Let $G\sq\dn\ti\dn.$
Geometrically, $G$ is a Borel set 
with countable vertical sections.
However, according to the Luzin--P.\,S.~Novikov 
theorem \cite{Kdst}, 
1) the projection of any Borel set 
with countable sections is itself a Borel set, 
so the domain $|G|$ 
of all vertices of the graph $G$ is equal 
to the projection 
$|G|=\ens{x}{\sus y\,(x\mG y)}$ of $G$ and 
is a Borel set, and 2) $G=\bigcup_nf_n$ is
a countable union of \rit{uniform} Borel sets 
$f_n$, that is, essentially, of graphs of Borel
functions, with the same domain  
(\ie, projection) 
$\dom f_n=|G|$. 
It remains to put $F=\ens{f_n}{n<\om}$.
\epf

As for countable Borel colorings, the next
lemma, along with Lemma~\ref{pro1}, show that
the implication of Lemma~\ref{22} for the 
Borel case
is irreversible.

\ble
\lam{pro2}
There are locally countable Borel graphs\/ $G$
with an uncountable Borel chromatic number. 
\ele

\bpf
Some examples are well known \eg\ from 
\cite{chrom,marks23} and others. 
An illustrative example can be given by 
the \rit{Vitali graph} 
$x\mV y$ iff $x\ne y$ and $x-y\in\dQ$,
defined on the real line $\dR$. 
Assuming that $c:\dR\to\om $ is a countable
Borel coloring for $V$, all sets 
$X_n=c\obr(n)$ are Borel (partial) 
transversals for the Vitali equivalence relation,
hence meager sets (as well as Lebesgue null sets). 
But this contradicts the Baire theorem
since $\dR=\bigcup_n X_n$.
\epf

The graphs of the second projective class $\fs12$
admit the following lemma  
about generating by a family of $\ali$ functions. 
Theorem~\ref{mts} shows that the 
generation by a \rit{countable} family of ROD  
functions (including all projective ones)
is, generally speaking, unprovable here.

\ble
\lam{t2}
Assume that\/ $G\sq\dn\ti\dn$ is a locally 
countable $\fs12$ graph. 
Then\/ $G$ is generated by a 
family\/ 
$F=\enx{f_\al}{\al<\omi}$ of\/ $\ali$ 
$\fd12$ functions.\snos
{A more thorough analysis of the G\"odel 
constructibility in the field of 
hereditarily countable sets makes it possible
to strengthen this result to the generation 
by a family  of $\ali$ 
\rit{Borel} functions, but this is beyond
the technical scope of this article.}
\ele
\bpf
Let $G$ be a $\is12(p_0)$ set; 
$p_0\in\dn$ is fixed in the proof.  
If $x\in\dn$ then the set 
$G(x)=\ens{y}{x\mathrel G y}$ 
of all $G$-adjacent elements is a countable 
$\is12(p_0,x)$ set, so that $G(x)\sq\rL[p_0,x]$. 
It follows that 
$G(x)\sq\ens{h(\al, p_0, x )}{\al<\omi}$, where 
$h(\al, p, x)$ is the $\al$th element of the set
$\dn\cap\rL[p,x]$ in the sense of the canonical 
G\"odel \weo\ of $\dn\cap\rL[p,x]$. 
We recall that  
$h:\omi\ti\dn\ti\dn\to\dn$ is known to be 
a $\id\hc1$ function. 

We let $f(\al, x)=h(\al, p_0, x)$ 
whenever $x\mG h(\al, p_0, x)$. 
Thus $f$ is function defined on 
$ 
D=\dom f :=
\enx{\ang{\al,x}\in\omi\ti\dn}
{x\mathrel G h(\al, p_0, x)},
$ 
and in fact a $\is\hc1(p_0)$ function, 
as the set of tuples 
$\ens{\ang{\al,x,f(\al, x)}}{\ang{\al,x}\in D}$, 
by the above.  
It follows that every $f_\al(x):=f(\al,x)$ 
is a $\fs12$ function from 
$D_\al=\ens{x\in\dn}{\ang{\al,x}\in D}$ to $\dn.$
Now prove that $G$ is generated by the  
family\/ 
$F=\enx{f_\al}{\al<\omi}$.

Assume that $x\mathrel G y$, hence $y\in G(x)$. 
Then $y=h(\al,p_0,x)=f(\al,x)=f_\al(x)$ for 
some $\al<\omi$ by construction. 
Conversely suppose that 
$y=f_\al(x)=f(\al,x)$ 
for some $\al$; then we have $x\mathrel G y$. 
\epf

We finish with a lemma related to the graph 
of equi-constructivity $\cL$. 
It is clear that, 
in the most constructive universe $\rL$, 
the graph $\cL$ is complete, \ie\ 
${x\ne y}\imp{x\mcL y}$. 
So $\cL$ is generated in $\rL$, for example, 
by an (uncountable) family of all continuous 
bijections $f:\dn\to\dn,$ 
which are not idempotent anywhere, 
\ie\ $f(x)\ne x$, $\kaz x$. 
However, there is also a countable generating 
family in $\rL$,
as the next lemma shows. 
This result can be compared with 
Theorem~\ref{mts}, saying that the graph $\cL$ 
is not 
generated by a countable family of ROD functions 
in the Solovay model.

\ble
\lam{t2+}
It is true in\/ $\rL$ that the 
equi-constructibility\/ $\is12$ graph\/ $\cL$ 
is generated by a countable family\/ 
$F=\enx{f_n}{n<\om}$ of\/ 
$\id12$ functions, and it has the chromatic 
number\/ $\chi(\cL)=\gc=\ali$. 
\ele

\bpf
We argue in $\rL$. 
Let $\lc$ be the G\"odel \weo\ of $\rL$.  
If $x\in\dn$ then let 
$\ens{h(n,x)}{n<\om}$ be the $\lc$-least 
enumeration of the set $\ens{y\in\dn}{y\lc x}$. 
We put $f_n(x):=h(n,x)$. 
Separately for the $\lc$-least real $x_0\in\dn$ 
put $f_n(x_0):=x_1$, where $x_1\in\dn$ is the 
$\lc$-next element. 
\epf

\parf{The equi-constructibility graph in the 
Sacks model} 
\las{saxm}

We consider the $\is12$ equi-constructibility 
graph $x\mcL y$ iff $x\ne y$ but $\rL[x]=\rL[y]$ 
on $\dn.$ 
Lemma~\ref{sax} below presents a rather 
difficult negative result for the Sacks-generic 
extensions, which 
will be used in the proof of Theorem~\ref{mts}.

Recall that the Sacks forcing $\saf z$ for a model 
of the form $\rL[z]$ 
consists of all perfect trees $T\in\rL[z]$, 
$T\sq\bse.$ 
It adjoins a 
\rit{Sacks-generic real} $a\in\dn$ to $\rL[z]$. 
In the proof of Lemma \ref{sax},  
we will use the following well-known
properties of the Sacks forcing, 
for which see, for example, \cite{gesh}
or \cite[Sec.\,15] {jechmill}.

\bpro
\lam{sx}
Assume that $\zo\in\dn.$  
Then$:$\vim
\ben
\renu
\vyk{
\itlb{sx1}%
\sloppy
the forcing\/ $\saf\zo$ preserves 
\/$\rL[\zo]$-uncountability, 
\ie,  
$\omi^{\rL[\zo]}=\omi^{\rL[\zo,a]};$ 
}%
\itlb{sx2}%
if\/ $a\in\dn$ is\/ $\saf\zo$-generic 
over\/ $\rL[\zo]$, and\/ 
$b\in\dn\cap\rL[\zo,a]$ then either\/ 
$b\in\rL[\zo]$ or\/ 
$a\in\rL[\zo,b]$ and then\/ $b$
also is\/ $\saf\zo$-generic over\/ 
$\rL[\zo]\,;$ 

\itlb{sx3}%
the forcing\/ $\saf\zo$ is homogeneous, \ie, if\/
$S,T\in\saf\zo$, then the cones\/  
$K_S=\ens{S'\in\saf\zo}{S'\sq S}$ and\/ $K_T$ 
are\/ $\sq$-isomorphic in\/ $\rL[\zo]\,;$

\itlb{sx4}%
as a standard consequence of \ref{sx3} 
in forcing theory,
if some\/ $S\in\saf\zo$
forces a closed formula\/ $\Phi$ with
parameters only from\/ $\rL[\zo]$, then any
other tree\/ $T\in\saf\zo$
also forces $\Phi$.
\hfill\qed
\een
\epro
 
\ble
\lam{sax}
Suppose that\/ $\zo\in\bai,$ and\/ $\ao\in\dn$ 
is a Sacks-generic real over\/ $\rL[\zo]$. 
Let, in\/ $\rL[\zo,\ao]$, 
$F=\sis{f_n}{n<\om}$ be an\/ $\OD(\ans\zo)$ 
sequence of functions\/ $f_n:\dn\to\dn$. 
Then, in\/ $\rL[\zo,\ao]$, $\cL$ 
is\/ {\bf not} generated by\/ $F.$ 
\ele

\bpf
If $x,y\in\bai$ then let $\sko xy\in\bai$ be 
the stepwise concatenation, that is, 
$(\sko xy)(2k)=x(k)$ 
and $(\sko xy)(2k+1)=y(k)$, $\kaz k$. 

Suppose to the contrary that, in $\rL[\zo,\ao],$ 
the family $F$ generates $\cL$. 

All statements about forcing below in the course 
of the proof are related only to forcing 
$\saf\zo$ over $\rL[\zo]$.
%
We consider the tree 
$$
T_0=\ens{s\in\bse}{\kaz n=2k<\lh s\,(s(2k)=\zo(k))}
\in\saf\zo\,,
$$
so that $[T_0]=\ens{\sko\zo y}{y\in\dn}$. 

Under our assumptions, there is a formula 
$\vpi(n,x,y)$ with parameters $\zo$ and 
some $\al_1,\dots,\al_m\in\Ord$ not explicitly 
indicated, such that it is true in $\rL[\zo,\ao]$ 
that $f_n=\ens{\ang{x,y}}{\vpi(n,x,y)}$ for 
each $n$.  
Let $\barf n$ be the shorthand for 
$\ens{\ang{x,y}}{\vpi(n,x,y)}$. 
Then the following sentence 
\bde
\item[$\Phi:=$] $\kaz n\,(\barf n:\dn\to\dn)$ and  
the family $\ens{\barf n}{n<\om}$ generates $\cL$ 
\ede
holds in $\rL[\zo,\ao]$ by the above. 
It follows that any condition in $\saf\zo$, 
in particular $T_0$, forces 
$\Phi$ over $\rL[\zo]$ by 
item \ref{sx4} of Proposition~\ref{sx}. 

From now on, {\ubf we argue in $\rL[\zo]$}.
Let $\una$ be the canonical $\saf\zo$-name for 
the principal Sacks-generic real in $\dn.$  
We claim that 
\ben
\Aenu
\itlb{saxA}%
{if\/ $T\in\saf\zo,\,T\sq T_0$, and $n<\om$, 
then there exist stronger conditions\/ 
$T',T''\in\saf\zo,\,T',T''\sq T$, and numbers\/ 
$k$ and\/ $j'\ne j''$ such that 
$T'$ forces $\barf n(\una)(k)=j'$ and 
$T''$ forces $\barf n(\una)(k)=j''$ 
over\/ $\rL[\zo]$.}
\setcounter{saved}{\value{enumi}}
\een

Indeed suppose towards the contrary that 
$T\in\saf\zo$ is a counterexample. 
Then there is a real $y\in\rL[\zo]\cap\dn$ 
such that in fact $T$ forces $\barf n(\una)=y$ 
over $\rL[\zo]$. 
We know that $\ao$ is a Sacks-generic real  
over $\rL[\zo]$. 
It follows that there is a Sacks-generic real 
$a\in [T]$ over $\rL[\zo]$ -- by item 
\ref{sx3} of Proposition \ref{sx}. 
Then we have $\barf n(a)=y$ by the genericity, 
and hence $a\mcL y$ holds in $\rL[\zo,a]$ 
because $T_0$ forces $\Phi$. 
 
However $y\in\rL[\zo]$ whereas $a\nin\rL[\zo]$ 
by the genericity. 
Therefore $\rL[a]=\rL[y]$ is definitely impossible. 
This contradiction completes the proof of \ref{saxA}. 

The following is an easy corollary of \ref{saxA}.
\ben
\setcounter{enumi}{\value{saved}}
\Aenu
\itlb{saxB}%
{if\/ $S,T\in\saf\zo,\,S,T\sq T_0$, and $n<\om$, 
then there exist stronger conditions\/ 
$S',T'\in\saf\zo,\,S'\sq S,\,T'\sq T$, 
and numbers\/ $k$ and\/ $j\ne \ell$ such that 
$S'$ forces $\barf n(\una)(k)=j$ and 
$T'$ forces $\barf n(\una)(k)=\ell$ 
over\/ $\rL[\zo]$.}
\setcounter{saved}{\value{enumi}}
\een

Now, {\ubf still arguing in $\rL[\zo]$}, we 
prove that  
\ben
\setcounter{enumi}{\value{saved}}
\Aenu
\itlb{saxC}%
In $\rL[\zo]$, there exists a 
system of trees $U_s$ and tuples $r_s$, 
where $s\in\bse,$ 
and tuples $\vt^n_s$, 
where $s\in\bse$ and $n<\lh s$, 
satisfying the following 
conditions \ref{sax1}--\ref{saxf}:
\setcounter{saved}{\value{enumi}}
\setcounter{enumi}{0}%
\vomh
%
\nenu
\itlb{sax1}%
$U_s\in\saf\zo$, $U_s\sq T_0$,  
$r_s\in U_s$, $\vt^n_s\in\bse$; 

\itlb{sax2}%
if\/ $s\in\bse$ and $n<\lh s$ then 
$\lh \vt^n_s\ge\lh s$; 

\itlb{sax3}%
if\/ $s\in\bse$ then $\lh r_s\ge \lh s$ and 
$r_s$ is $\sq$-comparable with each $t\in U_s$;


\itlb{sax4}%
if $s\su t$ belong to $\bse$ and $n<\lh s$ 
then $r_s\su r_t$, $U_t\sq U_s$, $\vt^n_s\su\vt^n_t$;

\itlb{sax5}%
if\/ $s\ne t\in\bse$ and $n<\lh s=\lh t$ then 
$r_{s}$ and $r_{t}$ are 
$\sq$-incomparable, and 
$\vt^n_{s}$, $\vt^n_{t}$ are 
$\sq$-incomparable as well; 

\itlb{saxf}%
if\/ $s\in\bse$ and $n<\lh s$ then 
$U_s$ forces $\vt^n_s\su \barf n(\kn a)$ 
over $\rL[\zo]$. 
\een
The construction {\ubf goes on in $\rL[\zo]$} 
by induction on $\lh s$. 

We put $U_\etu=T_0$ and $r_\etu=\etu$, 
where $\etu\in\bse$ is the empty tuple 
and the perfect tree $T_0\in\saf\pu$ 
was chosen above.

Now suppose that $m<\om$, and 
$U_s,r_s,\vt^n_s$ have been defined for all $s\in2^m$ 
(dyadic tuples of length $m$) and $n<m$, and 
satisfy conditions \ref{sax1}--\ref{saxf}.\vom 

{\sl Step 1.} 
Do the following for each $s\in2^m.$
Let $\rho_s=\roo{U_s}$ 
be the largest 
tuple $\rho\in U_s$ $\sq$-comparable with each 
$u\in U_s$. 
Then $r_s\sq\rho_s$ by \ref{sax3}. 
For $i=0,1$ define $r_{s\we i}=\rho_s\we i$ and 
$\baU_{s\we i}=
\ens{u\in U_s}{u\sq r_{s\we i}\lor r_{s\we i}\su u}$. 
Thus $r_\sg$, $\baU_\sg$ are defined for all 
$\sg\in2^{m+1}.$
The values $r_\sg$ are final, 
whereas the trees $\baU_\sg$ are temporary; 
they will be 
shrinked at the following steps. 
Note that the relevant parts of 
\ref{sax1},\ref{sax3},\ref{sax4},\ref{sax5} 
transfer to the level $m+1$ from level $m$, 
and will hold after any shrinking of the 
trees $\baU_\sg$ within $\saf\zo$.\vom

{\sl Step 2.} 
Do the following for all $s\in2^m$, $i=0,1$, $n\le m$. 

Put $\sg=s\we i$.  
Define a temporary tuple $\bavt^n_\sg$ as follows. 
Recall that $T_0$ forces $\Phi$, hence so does 
$\baU_\sg$, 
in particular $\baU_\sg$ forces 
$\barf n:\dn\to\dn.$ 
Therefore there is a tree 
$U\in\saf\zo$, $U\sq \baU_\sg$, 
and a tuple $\vt\in\bse$ 
with $\lh\vt>m$ and, if $n<m$ strictly 
(so that $\vt^n_s$ has been defined) 
then $\lh{\vt^n_s}<\lh\vt$, and in addition 
$U$ forces $\vt\su \barf n(\kn a)$.

We let $\bavt^n_\sg$ be such an $\vt$, 
and let the associated $U$ be 
the ``new'' $\baU_\sg$. 

This definition obviously honors \ref{sax2}, 
the last claim in \ref{sax4}, 
and \ref{saxf}.\vom

{\sl Step 3.} 
To fix \ref{sax5}, 
do the following for all $n\le m$ and  
$\sg\ne \ta$ in $2^{m+1}.$ 

Arguing as above (Step 2) and using \ref{saxB}, we 
find trees $U,U'\in\saf\zo$, $U\sq\baU_\sg$, 
$U'\sq\baU_\ta$, and tuples $\vt,\vt'\in\bse$ 
with $\bavt^n_\sg\su\vt$ and $\bavt^n_\ta\su\vt'$, 
such that $U$ forces $\vt\su \barf n(\kn a)$, 
$U'$ forces $\vt'\su \barf n(\kn a)$, 
and (this is where \ref{saxB} works!) 
for some $k<\lh\vt,\lh\vt'$ and $j\ne\ell$ we 
have $\sg(k)=j$ and $\sg'(k)=\ell$. 

The latter condition implies that 
tuples $\vt,\,\vt'$ 
are $\sq$-incomparable.

Let $\vt$ be the ``new'' $\bavt^n_\sg$, 
$\vt'$ be the ``new'' $\bavt^n_\ta$,
$U$ be the ``new'' $\baU_\sg$,
$U'$ be the ``new'' $\baU_\ta$. 
Go to the next triple of $n\le m$ and   
$\sg\ne \ta$ in $2^{m+1}.$\vom 

{\sl Step 3 -- finalization.} 
After processing all triples of $n\le m$ and  
$\sg\ne \ta$ in $2^{m+1},$ 
we let $U_\sg$ be the final tree $\baU_\sg$, and 
let $\vt^n_\sg$ be the final tuple $\bavt^n_\sg$ 
--- 
for all $n\le m$ and $\sg\in 2^{m+1}.$\vom 

{\sl Step 4 -- conclusion.} 
One easily sees that this construction yields 
a system of trees $U_s$ and tuples 
$r_s,\,\vt^n_s$ satisfying \ref{sax1}--\ref{saxf}. 
This ends the proof of \ref{saxC}.\vom

To make use of this system, 
{\ubf still arguing in $\rL[\zo]$}, we consider 
the tree  
\bce
$U=\ens{r\in\bse}
{\sus s\in\bse(r\sq r_s)}$; 
\ 
$U\sq U_\etu\sq T_0$ by construction.
\ece
Note that $U$ is a perfect tree by \ref{sax5} 
(regarding $r_s,r_t$), hence $U\in\saf\zo$. 
Let us prove two claims in connection 
with this tree $U$.

\ben
\setcounter{enumi}{\value{saved}}
\Aenu
\itlb{saxD}%
$\saf\zo$ forces, over $\rL[\zo]$,  
that ``each $\barf n$ 
is 1--1 on the set $X_U$ of all reals $x\in[U]$ 
$\saf\zo$-generic over $\rL[\zo]$''. \ 
\setcounter{saved}{\value{enumi}}
\een

Indeed, {\ubf arguing in a $\saf\zo$-generic 
extension $\rL[\zo,a]$ of $\rL[\zo]$}, consider 
any $n<\om$ and any reals $x\ne y$ in $X_U$. 
It follows from the definition of $U$ and 
\ref{sax3} that 
$[U]=\bigcap_m\bigcup_{s\in2^m}[U_s]$. 
Then, as $x\ne y$,  
there exist $m>n$ and $s\ne t$ in $2^m$ such that 
$x\in [U_s]$ and $y\in [U_t]$. 
Then the tuples $\vt^n_{s}$ and $\vt^n_{t}$ are 
$\sq$-incomparable by \ref{sax5}. 
On the other hand,  by \ref{saxf}, we have 
$\vt^n_s\su f_n(x)$ and $\vt^n_t\su f_n(y)$. 
We conclude that $f_n(x)\ne f_n(y)$, as required.  

\ben
\setcounter{enumi}{\value{saved}}
\Aenu
\itlb{saxE}%
$\saf\zo$ forces   
that ``(a) the set $X_U$ as in \ref{saxD} is 
uncountable, and 
(b) $X_U$ consists of pairwise
$\cL$-adjacent reals''. \ 
\setcounter{saved}{\value{enumi}}
\een

Indeed (a) follows from item \ref{sx2} of 
Proposition \ref{sx}, because it immediately 
implies $X_U=[U]\bez\rL[\zo]$. 

To prove (b) note that if $x,y$ are $\saf\zo$ generic 
over $\rL[\zo]$ then $\rL]\zo,x]=\rL[\zo,y]$ 
by item \ref{sx2} of Proposition \ref{sx}.
However, by construction $U\sq T_0$ and hence any 
real $z\in[U]$ satisfies $\zo\in\rL[z]$. 
Therefore the equality  $\rL]\zo,x]=\rL[\zo,y]$ 
implies $\rL[x]=\rL[y]$, hence $x\mcL y$, provided 
$x\ne y$ belong to $[U]$, as required. 

Thus  
it holds in the model $\rL[\zo,\ao]$ of 
Lemma~\ref{sax} by \ref{saxD},\ref{saxE}, 
that there is an uncountable set 
$X\sq\dn$ of 
$\mcL\msur$-adjacent elements, on which each $f_n$ 
is 1--1.  
Therefore, if $x_0\in X$ is fixed, then any  
$x\ne x_0$ in $X$ satisfies $x=f_n(x_0)$ or 
$x_0=f_n(x)$. 
Thus the set of all adjacent elements is countable 
due to the bijectivity of each $f_n$. 
This contradiction ends the proof 
of Lemma~\ref{sax}.
\epf

\parf{The equi-constructibility graph in the 
Solovay model} 
\las{solm}

Here we prove Theorem~\ref{mts}. 
The following definition introduces a particular 
form of the Solovay model we deal with in this 
theorem. 

\bdf
\lam{LN}
Let $\rL$ be the ground model and $\Om\in\rL$ 
be an inaccessible cardinal in $\rL$. 
Following \cite[Ch.\,26]{jechmill}, 
\cite[Ch.\,8]{Schi2014}, \cite{sol}, 
we let $\gN$ be the 
Levy--Solovay $\text{Coll}(\om,{<}\Om)$-generic 
extension of $\rL$; this is a model of $\ZFC$. 
\edf

The next proposition presents three well-known 
properties of the Solovay model which we'll 
use in the proof of Theorem~\ref{mts} below.

\bpro
[see {\cite{jechmill,Schi2014,sol,kl48}}]
\lam{sm}
The following is true in the model\/ $\gN$ 
of Definition~\ref{LN}$:$\vimh
\ben
\renu
\itlb{sm1}%
if sets\/ $X_0,X_1,X_2,\dots$ are 
real-ordinal definable\/ {\rm(ROD, for brevity)}, 
then the sequence\/ $\sis{X_k}{k<\om}$ is\/ 
$\ROD$ as well$;$ 

\vyk{
\itlb{sm2}%
if\/ $f:\dn\to\bai$ is a\/ 
$\ROD$ function then there 
is a perfect tree\/ $T\sq\bse$ such that\/ 
$f\res[T]$ is continuous$;$ 
}

\itlb{sm2}%
if\/ $z\in\dn$ and a set\/ $X\sq\rL[z]$ is\/ 
$\OD(z)$ 
{\rm(\ie, ordinal-definable 
with $z$ as an extra parameter)}
then\/ $X\in\rL[z]\,;$ 

\itlb{sm3}%
if\/ $a\in\dn$ and a countable set\/ $X\sq\bai$ 
is\/ $\OD(a)$ then\/ $X\in\rL[a]\,;$

\itlb{sm4}%
if\/ $z,a\in\dn$ 
and a set\/ $X\sq\bai$ is\/ 
$\OD(z)$ then the set\/ $X'=X\cap\rL[z,a]$ 
belongs to\/ $\rL[z,a]$ and is\/ 
$\OD(z)$ in\/ $\rL[z,a]\,;$

\itlb{sm5}%
if\/ $X\sq\dn$ is\/ $\ROD$ then either\/ 
$X$ is at most countable, or\/ $X$ has a perfect 
subset, and then\/ $X$ has the cardinality of  
continuum\/ $\gc=2^{\alo}\,;$

\itlb{sm6}%
if\/ $z\in\dn$ then the set\/ 
$\rL[z]\cap\bai$ is countable $($in\/ $\gN\,)$.
\hfill\qed
\een
\epro

\bpf[Theorem~\ref{mts}]
We skip the well-known parts of the theorem. 
For instance the local countability follows 
from Proposition~\ref{sm}\ref{sm6}.  

Let's focus on the key non-generation claim.
\rit{We argue in $\gN$}. 

Fix a family $F=\ens{f_n}{n<\om}$ of   
arbitrary ROD functions $f_n:\dn\to\dn.$ 
(If some $f_n$ originally has $\dom f_n=D\sneq\dn$ 
then we extend it by $f_n(x)=x^-$ for $x\nin D$, 
where $x^-(k)=1-x(k)$ for all $k$.) 
Suppose towards the contrary that $\cL$ 
is generated by $F.$ 
We observe that the whole sequence 
$S=\sis{f_n}{n<\om}$ is ROD (in $\gN$) by 
Proposition~\ref{sm}\ref{sm1}, hence  
there is a single real 
$z_0\in\dn$ such that $S$ is $\OD(z_0)$. 
Fix such a real $z_0$. 

\rit{Still arguing in $\gN$}, 
we also fix a real $\ao\in\dn$ Sacks-generic 
over $\rL[\zo]$.

Then by \ref{sm4} of Proposition~\ref{sm} each 
$f'_n=f_n\cap\rL[\zo,\ao]$ belongs to 
$\rL[\zo,\ao]$ and is $\OD(\zo)$ in 
$\rL[\zo,\ao]$, 
and moreover, the whole sequence 
$S'=\sis{f'_n}{n<\om}$ belongs to 
$\rL[\zo,\ao]$ and is $\OD(\zo)$ in 
$\rL[\zo,\ao]$. 

On the other hand, 
given any $x\in\rL[\zo,\ao]\cap\dn,$ 
we have  
$f_n(x)\in\rL[\zo,\ao]$ and 
$f_n\obr(x)\sq\rL[z_0,\ao]$ 
(because $f_n\obr(x)$ is countable) 
by resp.\ \ref{sm2} and \ref{sm3} of 
Proposition~\ref{sm}. 
It follows that 
$f'_n=f_n\res{(\dn\cap\rL[z_0,\ao])}$, 
and hence it is true in $\rL[z_0,\ao]$ 
that the $\OD(\zo)$ sequence $S'$ 
generates $\cL.$
But this contradicts Lemma~\ref{sax}\,! 
This completes the proof of the main part 
of Theorem~\ref{mts}. 

To prove the additional statement
of Theorem~\ref{mts}, it is sufficient 
to refer to Lemma~\ref{22}.
\epf

\vyk{
Now prove the ``moreover'' part.
We have to show that\/ $\cL$ has  
the\/ $\ROD$ chromatic number equal to\/ 
$\gc=2^{\alo}$ in the Solovay model\/ $\gN$.

Indeed, 
otherwise the\/ $\ROD$ chromatic number of $\cL$ 
is ${\le}\alo$ by Proposition~\ref{sm}\ref{sm5}, 
meaning that, in $\gN$, there is a $\ROD$ map 
$c:\dn\to\om$ satisfying 
(*) ${x\cL y}\imp c(x)\ne c(y)$. 
\rit{Arguing in $\gN$}, we let 
$U_x=\ens{c(y)}{y\in[x]_\rL}$ for any real 
$x\in\dn,$ where 
$[x]_\rL=\ens{y\in\dn}{\rL[x]=\rL[y]}$, the 
\rit{constructibility degree} of $x$. 
Thus $U_x\sq\om$ is infinite by (*) 
for each $x$. 

\rit{Still arguing in $\gN$}, define $f(x)\in\dn$ 
for any $x\in\dn$ as follows. 
Let $n=c(x)$. 
As $U_x$ is infinite, let $n'>n$ be the immediate 
successor of $n$ in $U_x$. 
There is a unique $y\in[x]_\rL$ with $c(y)=n'$. 
Put $f(x)=y$. 
Then put 
$f_k(x)=
\underbrace{f(f(\dots}_{k\text{ times }f}(x)\dots))$.
Then $f_k:\dn\to\dn$ are $\ROD$ maps, and easily 
$\cL$ is generated by $F=\ens{f_k}{k\ge1}$. 
But this contradicts the main part 
of Theorem~\ref{mts}.  
}

\parf{The case of $\fp12$ graphs} 
\las{p12}

Our {\ubf proof of Theorem~\ref{mt}} here 
is based 
on a model defined in \cite{kl28}, in which 
there exists a non-$\ROD$-uniformizable $\ip12$ 
planar set with countable cross-sections.
For the convenience of the reader, 
we present here this construction, without going 
into technical details, and then show how to 
convert it into an example for Theorem~\ref{mt} 
in the same model. 

Beginning with $\rL$ as the ground model, we 
defined in \cite[\S\,9]{kl28} a sequence 
$\sis{\dP_\xi}{\xi<\omi}\in\rL$ of forcing 
notions $\dP_\xi$. 
Each of $\dP_\xi$ consists of perfect trees 
$T\sq 2^{<\om}$ and is rather 
similar to the Jensen minimal-$\id13$-real 
forcing considered in detail \eg\ in 
\cite[\S\,28.A]{jechmill} or \cite{jml19}. 

Then the finite-support product 
$\dP=\prod_{\xi<\omi}\prod_{k<\om}\dP_{\xi k}
\in\rL$ is defined in \cite{kl28}, 
where each factor 
$\dP_{\xi k}$ is equal to $\dP_\xi$, 
and we proved the following there:  

\ben
\nenu
\itlb{e3}\msur%
$\dP$ does not collapse $\rL$-cardinals; 

\itlb{e4}\msur%
$\dP$ adjoins a generic array 
$X=\sis{x_{\xi k}}{\xi<\omi,\,k<\om}$
of reals $x_{\xi k}\in 2^\om$; 

\itlb{e5}%
each $x_{\xi k}$ is $\dP_\xi$-generic over 
$\rL$, and conversely, every real $x\in\rL[X]$,   
$\dP_\xi$-generic over $\rL$, is equal to one 
of $x_{\xi k},\:k<\om$; 

\itlb{e6}%
the relation \ 
$\text{\rm Gen}(\xi,x):=$ 
``$\xi<\omi$ and 
$x\in2^\om$ is $\dP_\xi$-generic over $\rL$'' \ 
is $\ip\HC1$ in $\rL[X]$, where $\HC=$  
all hereditarily countable sets$;$ 

\itlb{e7}%
by the finite-support product forcing 
theory, if $A\in\rL$, $A\sq\omi\ti\om$, and 
$\ang{\xi,k}\nin A$ then 
$x_{\xi k}\nin \rL[X\res A]$, where 
$X\res A=\sis{x_{\xi k}}{\ang{\xi,k}\in A}$, 
and moreover, $x_{\xi k}\nin \OD(X\res A)$ 
in $\rL[X]$. 
\setcounter{saved}{\value{enumi}}
\een
We used these properties of the generic model 
$\rL[X]$ in 
\cite{kl28} to prove that the set 
$W=\ens{\ang{\xi,x_{\xi k}}}{\xi<\omi\land k<\om}$
\vyk{
\setcounter{equation}{\value{saved}}
\bul
{e8}
{W=\ens{\ang{\xi,x_{\xi k}}}{\xi<\omi\land k<\om}.
\setcounter{saved}{\value{equation}}}
}
%
is a non-$\ROD$-uniformizable 
$\ip\HC1$ set with countable sections 
$W_\xi=\ens{x_{\xi k}}{k<\om}$ in $\rL[X]$. 

This set was easily converted in \cite{kl28} to a 
$\ip12$ set $W'\sq\dn\ti\dn$ in $\rL[X]$ 
with the same properties. 
Indeed let $\text{\ubf WO}\sq\bai$ be the standard 
$\ip11$ set of codes for countable ordinals, and 
if $w\in \text{\ubf WO}$ then let $|w|<\omi$ be 
coded by $w$. 
The set 
$W'=\ens{\ang{w,x_{\xi  k}}}
{w\in \text{\ubf WO}\land |w|=\xi\land k<\om}$ 
is then a non-$\ROD$-uniformizable 
$\ip12$ set with countable sections in $\rL[X]$. 

A similar design of forcing, but in the form
of a simpler product 
$\prod_{\xi<\omi}\dP_{\xi}$,
was used in \cite{kl30} for other purposes.

\vyk{
If $p,q\in \dn\ti\dn$ then define $p\rE q$ iff 
either $p=q$, or both  
$p=\ang{w,x_{\xi  k}}$ and 
$q=\ang{w',x_{\xi', n}}$ belong to $W'$ 
and $w=w'$ (then $\xi=\xi'$), 
but not necessarily $k=n$. 
Then, in $\rL[X]$, 
$\rE$ is a $\ip12$ equivalence relation 
on $\dn$ (as $W'$ is $\ip12$)  
with countable equivalence classes 
since 
$[\ang{w,x_{\xi  k}}]_{\rE}=
\ens{\ang{w,x_{\xi  j}}}{j<\om}$. 
}%

\bpf[Theorem~\ref{mt}] 
If $p=\ang{w,x_{\xi k}}$ and 
$q=\ang{w',x_{\xi', n}}$ belong to $W'$ 
then define $p \mG q$ iff 
$w=w'$, $\xi=\xi'$, and $k\ne n$. 
Thus, in $\rL[X]$,  
$G$ is a locally countable graph 
of class $\ip12$ 
with $W'=|G|$ as the set of vertices.

To complete the proof of Theorem~\ref{mt}, 
it remains to show that, in $\rL[X]$, $G$ is 
not generated by a countable family 
$F=\ens{f_n}{n<\om}$ of ROD functions 
$f_n:W'\to W'$. 
Assume to the contrary that $G$ is generated 
by such an $F.$ 

\rit{Arguing in\/ $\rL[X]$}, for any $n$ there 
is a real $u_n\in\dn,$ such that $f_n$ is 
$\OD(\ans{u_n})$. 
Then by \ref{e3} above there also exists a 
countable set $A_n\sq \omi\ti\om$ with 
$u_n\in\rL[X\res A_n]$. 
The set $A=\bigcup_nA_n$ is countable as well, 
hence there exist pairs 
$\ang{\xi,k}$ and 
$\ang{\xi,n}$ in $(\omi\ti\om)\bez A$, 
with the same $\xi$ and with $k\ne n$. 
Pick a code $w\in\text{\bf WO}$ with $|w|=\xi$.
Then the according elements 
$p=\ang{w,x_{\xi k}}$ and 
$q=\ang{w,x_{\xi n}}$ in $W'$ satisfy 
$p=f_m(q)$ or $q=f_m(p)$ for some $m$ 
by the contrary assumption above. 

Let say $p=f_m(q)$.
Then $x_{\xi k}$ is 
$\OD(\ans{w,x_{\xi n},u_m})$ in $\rL[X]$ 
because $f_m$ is $\OD(\ans{u_m})$. 
It follows that $x_{\xi  k}$ is 
$\OD(\ans{x_{\xi  n},X\res A})$ since $w\in\rL$. 
But this contradicts \ref{e7} as  
$\ang{\xi,k}\notin A\cup\ans{\ang{\xi,n}}$. 
This completes the proof of the main part 
of Theorem~\ref{mt}. 

To prove the additional statement
of Theorem~\ref{mt}, it is sufficient 
to refer to Lemma~\ref{22}.
\epf

\vyk{
As for the ``moreover'' part of the theorem, 
to deduce the uncountability of the 
$\ROD$-chromatic number of $G$, it suffices to 
repeat the related reasoning from the 
last part of the proof
of Theorem~\ref{mts}, changing $[x]_\rL$ 
to the connected $G$-component, defined for  
any $p=\ang{w,x_{\xi k}}\in W'$ so that:
$$
[p]_G=
\enx{\ang{w',x_{\xi', n}}\in W'}
{w=w'\land \xi=\xi'\land k,n<\om}.\eqno\qed
$$ 
}

\parf{Concluding remarks} 
\las{zz}

Our Theorem~\ref{mt} solves, in the positive, 
a problem on the existence of locally 
countable $\ip12$ graphs non-generated 
by countable 
families of definable functions. 
Two separate results are obtained for $\fs12$ 
graphs in \S\,\ref{s12}. 
Two non-generation results on the $\is12$ 
equi-constructibility graph $\cL$  
(Lemma~\ref{sax} and Theorem~\ref{mts})
are obtained in \S\S\,\ref{saxm},\ref{solm}. 
We expect that the results obtained and 
methods developed will find further 
applications in modern research in
descriptive set theory and forcing.

\vyk{
We also expect that our methods of effective
transfinite constructions will make a definite
contribution to the modern theory of generalized
computability on uncountable
structures and the theory of information 
transmission between
structures on ordinals and Borel
structures associated with the real line, 
as in recent works \cite{ham_effmat,hamTM}.
}


We finish with  
the following problems that arise from our study. 

\bvo
Is there a locally countable $\fp11$ 
graph in the Solovay model $\gN$ 
(as in Definition~\ref{LN}) 
not generated by a 
countable family of\/ $\ROD$ functions?
\evo

A possible plan to solve the problem could 
be as follows. 
Using the Novikov-Kondo-Addison 
uniformization, we can start with a $\ip11$ 
set $U\sq(\dn)^3$, uniform in the sense 
$(\dn\ti\dn)\ti\dn,$ and such that 
${x\mcL y}\eqv \sus p\,U(x,y,p)$. 
Consider the $\ip11$ graph $\cH$ whose domain 
satisfies $|\cH|\sq\dn\ti\dn,$ defined so that 
$\ang{x,p} \cH \ang{y,q}$ iff $p=q$ and 
$U(x,y,p)$. 
Clearly $\cH$ is locally countable in the 
Solovay model $\gN$. 
It remains to show that $\cH$ is not 
generated by a 
countable family of\/ $\ROD$ functions 
in $\gN$.

\bvo
Does Lemma~\ref{sax} remain true for reals 
$\ao$ Cohen-generic or Solovay-random 
(or any other popular type of generic reals)? 
\evo

\bvo
Find a simpler proof of Theorem 
\ref{mts}, circumventing Lemma~\ref{sax}.
\evo

The following argument may be suggested. 
\rit{Suppose to the contrary}, that the graph $\cL$
is generated by a family $F=\ens{f_n}{n<\om}$ 
of ROD functions $f_n:\dn\to\dn$ in the 
Solovay  model $\gN$. 
According to the general properties of the model, 
there is a perfect tree $T_0\sq\bse$ such that all 
restricted functions $f_n\res[T_0]$ are continuous. 
Further, according to the properties of the Cantor 
discountinuum $\dn,$
there is a perfect tree $T\sq T_0$ such  
that each $f_n\res[T]$ is either 1--1 or
a constant. 
The case of a constant is rejected by the local 
countability of the graph $\cL$ in the Solovay 
model $\gN$, and therefore
we assume that all $f_n\res[T]$ are 1--1.

There is a real $\zo\in\dn,$ such that both  
$T$ and a suitable sequence of codes
for continuous functions $f_n\res[T]$ 
belong to $\rL[\zo]$. 
Now, to deduce the contradiction as in the proof
of Theorem~\ref{mts} above, it would be 
enough to find an 
$\rL[\zo]$-uncountable set $X\sq[T]$ 
of $\cL$-adjacent 
elements, in $\rL[\zo]$ or in an  
$\rL[\zo]$-uncountability preserving generic 
extension
of the model $\rL[\zo]$. 
This seems doable in the case where $\zo$
preserves $\rL$-uncountability, 
but in the general case
it remains an interesting open question.

\bvo
Coming back to inequality \eqref{chrn} 
in \S\,\ref{prel}, we may ask whether there is a 
model of $\zfc$ in which every inequality there 
is strict for $\cL$ or any other graph.
\evo

\bvo
In view of Lemma~\ref{t2+} and 
Theorem~\ref{mts}, we may ask whether there is a 
model of $\zfc$ in which 
$\chi_{\text{\sc rod}}(\cL)=\alo$.
\evo

For $\chi_{\text{\sc bor}}$, a similar question
is likely to have a negative answer.

\vyk{
\begin{ack}
The authors are grateful to Mirna D\v zamonja for 
an interesting
discussion and valuable comments.
\end{ack}
}

\let\section\subsection

\renek{\refname} {References}

\bibliographystyle{plain}

\small

\bibliography
{../0bib/all.bib,../0bib/kl.bib,../0bib/mar.bib}

\end{document}